\documentclass[11pt,a4paper]{article}
\usepackage[english]{babel}
\usepackage{amsmath,amsthm}
\usepackage{amsfonts}
\usepackage{amssymb}

\usepackage{hyperref}
\usepackage{a4wide}

\usepackage{tikz}
\usepackage{caption}

\def\doi{\begingroup\catcode`\_=12\relax\doiarg}
\def\doiarg#1{\href{https://doi.org/#1}{\texttt{doi:#1}}\endgroup}

\newtheorem{thm}{Theorem}[section]
\newtheorem{cor}[thm]{Corollary}
\newtheorem{lem}[thm]{Lemma}

\theoremstyle{definition}

\theoremstyle{remark}

\numberwithin{equation}{section}

\title{2-Factors in Graphs\,\thanks{\,The work leading to this paper was carried out at LSE, London, in 2011.
		Discussions with Monika M\o bjerg Andersen at the University of Southern Denmark in the spring of 2011 are also acknowledged.
		The second author thanks the Danish Council for Independent Research (FNU) and the University of Southern Denmark for support, and the Department of Mathematics at LSE for hospitality and a stimulating environment.}}
\author{Jan van den Heuvel\\[1mm]Department of Mathematics,\\
	London School of Economics \& Political Science\\
	Houghton Street, London WC2A 2AE, UK\\
	\href{mailto:j.van-den-heuvel@lse.ac.uk}{\tt j.van-den-heuvel@lse.ac.uk}\\[2mm]
	and\\[2mm]
	Bjarne Toft\\[1mm]
	Department of Mathematics and Computer Science,\\
	University of Southern Denmark\\
	Campusvej 55, DK-5230 Odense M, Denmark\\
	\href{mailto:btoft@imada.sdu.dk}{btoft@imada.sdu.dk}}

\date{}

\begin{document}

\maketitle

\begin{abstract}
\noindent
An account of $2$-factors in graphs and their history is presented.
We give a direct graph-theoretic proof of the $2$-Factor Theorem and a new variant of it, and also a new complete characterisation of the maximal graphs without $2$-factors.
This is based on the works of Tibor Gallai on $1$-factors and of Hans-Boris Belck on $k$-factors, both published in 1950 and independently containing the theory of alternating chains.
We also present an easy proof that a $(2k+1)$-regular graph with at most $2k$ leaves has a $2$-factor, and we describe all connected $(2k+1)$-regular graphs with exactly $2k+1$ leaves without a $2$-factor.
This generalises Julius Petersen's theorem that any $3$-regular graph with at most two leaves has a $1$-factor, and it generalises the extremal graphs Sylvester discovered for that theorem.
\end{abstract}

\section{Introduction}

We use standard graph-theoretic notation, as can be found e.g.\ in Bondy and Murty~\cite{B&M 2008}, unless otherwise indicated.
Our graphs are undirected and finite, but we do allow multiple edges and loops.

Given a graph $G$ and $X,Y\subseteq V(G)$, we use $E(X,Y)$ to denote the set of edges with one end in $X$ and the other end in $Y$ in $G$, and $e(X,Y)$ to denote their number.
The subgraph of $G$ induced by $X$ is denoted $G[X]$.

An edge in a graph not contained in a cycle is a \emph{bridge} of the graph.
The removal of a bridge from a graph results in two new connected components replacing the connected component containing the bridge.
A \emph{leaf} of a graph is an induced bridgeless connected subgraph joined to the rest of the graph by one bridge.
Two different leaves of a graph are disjoint.
The number of leaves in a connected graph is at most the number of bridges, except that a connected graph with exactly one bridge has two leaves.

The degree of a vertex $x$ is denoted $d(x)$, where a loop adds $2$ to the degree of its vertex.
A \emph{$k$-factor} of a graph is a spanning subgraph of that graph in which all vertices have degree~$k$.

\medskip
Graph theory in the form we know it was created in the collaboration between James Joseph Sylvester in Oxford and Julius Petersen in Copenhagen.
The collaboration started with Sylvester's visit to Copenhagen in September 1889, and was followed by an exchange of letters and a visit by Petersen to England in December 1889 and January 1890, resulting in Petersen's famous paper \cite{Petersen 1891}.
The main results of the paper are the following two theorems

\begin{thm}[Petersen, 1891 \cite{Petersen 1891}]\mbox{}\\*
	Let $G$ be a $3$-regular graph with at most two leaves.
	Then $G$ has a $1$-factor.
\label{JP3}
\end{thm}

\begin{thm}[Petersen, 1891 \cite{Petersen 1891}]\mbox{}\\*
	Let $G$ be a $2r$-regular graph, for some positive integer $r$.
	Then $G$ has a $2$-factor.
\label{JP2}
\end{thm}

\noindent
It is often mentioned that graph theory originated with Euler's solution of the K\"onigsberg Bridge Problem in 1736~\cite{Euler 1736}.
Euler observed (with a different terminology) that a connected graph has a closed walk without repeated edges containing all edges of the graph if and only if all the degrees of the vertices of the graph are even.
However, Euler did not define graphs, nor did he draw any graph.
Petersen did not read much and was probably not familiar with Euler's paper.
It seems that he obtained Euler's Theorem independently, and he observed that the case $r=2$ of Theorem~\ref{JP2} immediately follows by colouring the edges of an Euler walk alternating red and blue.

D\'enes K\H{o}nig wrote the first extensive monograph of graph theory in 1936 \cite{Konig 1936} and gave Petersen's results a prominent place.
K\H{o}nig presented different ways to prove and relate factorisation theorems for general and for bipartite graphs.
In particular, Theorem \ref{JP2} above and K\H{o}nig's own Theorem \ref{KonigBip} on bipartite graphs below may be considered equivalent since each may be deduced from the other.

\begin{thm}[K\H{o}nig, 1936 \cite{Konig 1936}]\mbox{}\\*
	Let $G$ be a regular bipartite graph.
	Then $G$ has a $1$-factor.
\label{KonigBip}
\end{thm}

\noindent
In 1947, William T. Tutte obtained the general $1$-Factor Theorem.

\begin{thm}[{\rm\textbf{1-Factor Theorem}}; Tutte, 1947 \cite{Tutte 1947}]\mbox{}\\*
	A graph $G$ has a $1$-factor if and only if for every $A \subseteq V(G)$ the graph $G-A$ has at most $|A|$ connected components with an odd number of vertices.
\label{1FT}
\end{thm}

\noindent
Tutte's proof is based on a characterisation of the maximal graphs $G$ without a $1$-factor (called \emph{hyperprime graphs} in Tutte's paper \cite{Tutte 1947}).
If $G$ is maximal without a $1$-factor, $|V(G)|$ is even, and  $B$ denotes the set of vertices joined to all other vertices, then the connected components of $G-B$ are all complete.
Theorem \ref{1FT} follows easily from this characterisation.
Tutte's proof of this characterisation is algebraic, using skew-symmetric determinants.
In 1950 Hans-Boris Belck~\cite{Belck 1950} and Tibor Gallai~\cite{Gallai 1950} independently published graph theoretic proofs, while in 1952 Frederic George Maunsell~\cite{Maunsell 1952} replaced the algebraic part of Tutte's proof by arguments using graph properties only.
In 1975 L\'aszl\'o Lov\'asz \cite{Lovasz 1975} gave a short graph-theoretic proof.

Gallai \cite{Gallai 1950} obtained a proof of Theorem \ref{1FT} based on his general theory of alternating chains.
Gallai's  method of proof gives a structural description of all graphs in terms of their maximum matchings, now known as the Gallai-Edmonds Decomposition Theorem (see \cite{Lovasz 1986}).
Applying the theory of alternating chains, Gallai generalised  in elegant fashion many of the known factorisation theorems and proved new ones.
Among other results, he generalised  Baebler's theorem \cite{Baebler 1938} that any $(2r+1)$-regular graph without bridges contains a $2$-factor.
Gallai's paper has an intriguing footnote, that he had obtained a similar structural theory for general graphs (as opposed to regular graphs) in terms of their $2$-factors: ``\emph{With the present method I have succeeded in getting factorisation theorems for general graphs besides $\sigma =1$ only for the case $\sigma =2$. I shall discuss these results on another occasion.}''

Unfortunately, Gallai never published his theory of $2$-factors.
One may guess that the reason was that he discovered that Belck \cite{Belck 1950}, also in 1950, published a general $k$-factor theorem, generalizing Theorem \ref{1FT}.
As Belck's paper shows, and we explain in this paper, the theory of alternating chains is the key to a proof of the general $k$-factor theorem.
For a graph without a $k$-factor, Gallai considered an induced subgraph with all degrees at most $k$ with minimum deficiency from being $k$-regular, whereas Belck added new edges to get a maximal graph without a $k$-factor (called \emph{hyper-$k$-prime graphs} in \cite{Belck 1950}), like Tutte did for $1$-factors in 1947.
Thus Belck and Gallai independently gave graph-theoretic proofs of Tutte's $1$-Factor Theorem, and Belck was the first to obtain the general $k$-factor theorem.

\begin{thm}[Belck, 1950 \cite{Belck 1950}]\mbox{}\\*
	Let $k$ be a positive integer.
	A graph $G$ has a $k$-factor if and only if for every disjoint $A,B\subseteq V(G)$, setting $C=V(G)\setminus(A\cup B)$, the graph $G[C]$ has at most $2e(A,A)-k|A|+ k|B|+e(A,C)$ connected components $C_i$ with $k|V(C_i)|+e(A,V(C_i))$ odd.
\label{kFT}
\end{thm}

\noindent
In particular, for $2$-factors we obtain the following.

\begin{cor}\mbox{}\\*
	A graph $G$ has a $2$-factor if and only if for every disjoint $A,B\subseteq V(G)$, setting $C=V(G)\setminus(A\cup B)$, the graph $G[C]$ has at most $2e(A,A)-2|A|+ 2|B|+e(A,C)$ connected components joined to $A$ by an odd number of edges.
\label{2FT}
\end{cor}

\noindent
As we shall see, in the condition of this corollary one may consider only independent sets $A$.

\begin{thm}[\rm\textbf{2-Factor Theorem}]\mbox{}\\*
	A graph $G$ has a $2$-factor if and only if for every disjoint $A,B\subseteq V(G)$, $A$ an independent set, and setting $C=V(G)\setminus(A\cup B)$, the graph $G[C]$ has at most $-2|A|+ 2|B|+e(A,C)$ connected components joined to $A$ by an odd number of edges.
\label{2FTrevised}
\end{thm}

\noindent
In 1952, Tutte \cite{Tutte 1952}, citing the works of Belck \cite{Belck 1950} and Gallai \cite{Gallai 1950}, subsequently obtained a general $f$-factor theorem, giving a necessary and sufficient condition for the existence of an $f$-factor in a general graph (where $f$ is a mapping of $V(G)$ into the non-negative integers).
A simple proof of the $f$-factor theorem, deducing it from the $1$-Factor Theorem \ref{1FT}, is also due to Tutte~\cite{Tutte 1954}.
This attractive proof may also be found in the book by Bollob\'as \cite{Bollobas 1978}.
Already Claude Berge in 1958, in `the second book' on graph theory \cite{Berge 1958}, described Gallai's theory and used Tutte's method to prove the $2$-Factor Theorem (formulated in \cite{Berge 1958} with misprints), calling $2$-factors \emph{semi-factors}.
New comprehensive texts on factorisation, containing the general \mbox{$f$-factor} theorem, are due to Yu and Liu \cite{Yu 2009} and to Akiyama and Kano \cite{Akiyama 2011}.

For maximum matchings in general graphs, Jack Edmonds \cite{Edmonds 1965} in 1965 obtained a polynomial algorithm, showing that the ``good'' characterisation of Tutte \cite{Tutte 1947} is accompanied by a ``good'' algorithm.
The output of the algorithm, with a general graph as input, gives, in a different language, the same structure as  Gallai's method~\cite{Gallai 1950}.
That structure is now usually called the \emph{Gallai-Edmonds decomposition} of a graph.

\subsection{Our work}

Since an edge-multiplicity of $3$ or more, or two or more loops at the same vertex will not help to create a $2$-factor, from now on we will only consider the class $\mathcal{M}_2$ of graphs having edge-multiplicity at most $2$ and each vertex having at most one loop.

One aim of the present paper is to describe the complete structure of the maximal graphs without $2$-factors.
None of the above cited works contains such a complete description.
Since cycles in graphs and $2$-factors are basic structures in graph theory, it seems appropriate to fill this gap.

\begin{thm}\mbox{}\\*
	Let $\mathcal{M}_2$ be the class of all graphs with all multiplicities at most $2$ and with each vertex having at most one loop.
	Then $G$ is maximal without a $2$-factor within the class $\mathcal{M}_2$ if and only if the following is satisfied.

	Let the set~$A$ contain all vertices of $G$ without loops, the set $B$ contain all vertices in $G$ with loops and joined to every other vertex in $G$ by two edges, and set $C=V(G)\setminus(A\cup B)$.
	Suppose $G[C]$ has~$q$ connected components $\mathcal{C}=\{C_1,\ldots,C_q\}$.
	Then we have the following.
	\vspace{-0.6mm}
	\begin{itemize}
		\addtolength{\itemsep}{-2mm}
		\item The set $A$ is independent.
		\item The components in $\mathcal{C}$ are all complete (with a loop at every vertex and two edges between any two vertices in the same component).
		\item Each component $C_i$ is joined to $A$ by an odd matching (of size $e(A,V(C_i))$).
		\item $2|A|+q=2|B|+2+e(A,C)$, which is equivalent to $|A|= |B|+1+\sum\limits_{i=1}^q\tfrac12\bigl(e(A,V(C_i))-1)$.
		\item For all $A'\subseteq A$, $A'\ne\varnothing$, and all $\mathcal{C}'\subseteq\mathcal{C}$ we have: $2|A'|+|\mathcal{C}'|\ge 2+\sum\limits_{C_i\in \mathcal{C}'}e(A',V(C_i))$.
	\end{itemize}
	\label{final}
\end{thm}

\noindent
On the way to proving this theorem we shall give a simple direct proof of the  $2$-Factor Theorem, using the method of Belck and Gallai.
Moreover, we deduce in a simple way from the $2$-Factor Theorem the following result.

\begin{thm}\mbox{}\\*
	Let $G$ be a $(2k+1)$-regular graph with at most $2k$ leaves.
	Then $G$ has a $2$-factor.
\label{regular}
\end{thm}

\noindent
The case $k=1$ is Petersen's classical Theorem \ref{JP3} on $3$-regular graphs (since a $1$-factor is the complement of a $2$-factor in a $3$-regular graph).
Being a special case of a general $2$-factor theorem, Petersen's theorem should be regarded as a $2$-factor theorem, and not just a $1$-factor theorem.
According to Roland H\"aggkvist \cite{Roland 2001}, he announced the result in Theorem~\ref{regular} at a meeting in Oberwolfach in 1977.

Theorem~\ref{regular} with the stronger condition that there are no bridges is due to Baebler~\cite{Baebler 1938}, and in generalised forms (but still for graphs without bridges) to Belck \cite{Belck 1950} and Gallai~\cite{Gallai 1950}.
Hanson, Loten and Toft \cite{Hanson 1998} proved a result similar to Theorem~\ref{regular}, allowing up to $2k$ bridges.
In hindsight, their proof looks unnecessarily complicated now.
A better proof was obtained by Douglas West (private communication), and a generalisation was published in 2021 by Kostochka, Raspaud, Toft, West, and Zirlin \cite{KRTWZ 2021}.
That paper also contained a description of all $(2k+1)$-regular graphs with exactly $2k+1$ bridges without a $2$-factor; we also describe these graphs, but in a slightly different and generalised way.
These extremal graphs generalise the Sylvester graphs that Sylvester communicated to Petersen in 1889 \cite{Petersen 1891}.

The final section of our paper contains some of the information we discovered about Hans-Boris Belck, whose work seems to be often overlooked by those working on factors in graphs.
A reason might be that the paper was written in quite complicated German.

\section{The proof of the 2-Factor Theorem}

\subsection{The `only if' part of Corollary~\ref{2FT} and the 2-Factor Theorem}

Recall that $\mathcal{M}_2$ denotes the class of graphs having edge-multiplicity at most $2$ and at most one loop at each vertex.

Let $G$ be a graph from~$\mathcal{M}_2$ containing a $2$-factor~$F$.
Let $A$ and $B$ be two disjoint subsets of $V(G)$.
Set $C=V(G)\setminus(A\cup B)$ and let $C_1,C_2,\ldots$ be the connected components of $G[C]$, the subgraph of $G$ induced by $C$.

The factor $F$ consists of disjoint cycles.
Each such cycle may contain vertices from $A$, and between these vertices there are paths with vertices from $B\cup C$.
We shall distinguish four different types of these paths, and use these to construct three disjoint subsets $B_1,B_2,B_3$ of~$B$.
\vspace{-2mm}
\begin{itemize}
	\addtolength{\itemsep}{-2mm}
	\item Paths in which both end-vertices are from $B$ (which includes the case of a path with one vertex (from $B$)).
	For those paths, add one of the end-vertices to $B_1$.
	\item Paths in which both end-vertices are from the same component $C_i$ of $G[C]$.
	\item Paths in which the end-vertices are from two different components $C_i$ and $C_j$ of $G[C]$.
	Then the path must contain a vertex from $B$; add one such vertex from $B$ to $B_2$.
	\item Paths with one end-vertex from $B$ and one from $C$;
	add the end-vertex from~$B$ to~$B_3$.
\end{itemize}

\noindent
Note that the number of edges from $E(A,B)$ in $F$ is equal to $2|B_1|+|B_3|$.
Hence the number of edges from $E(A,B)$ not in $F$ is equal to $e(A,B)-2|B_1|-|B_3|$.

Let $q$ denote the number of components $C_i$ in $G[C]$ with an odd number of edges to $A$, i.e.\ with $e(A,V(C_i))$ odd.
At most $2|B_2|+|B_3|$ of these $q$ components can have all the edges to~$A$ contained in $F$, since the $(A,V(C_i))$-edges from the second type of paths above come in pairs.
Hence at least $q-2|B_2|-|B_3|$ of the $(A,C)$-edges are not in $F$.

From this we deduce
\begin{align*}
	&\sum_{x\in A}(d(x)-2)\\[-3mm]
	&\qquad\qquad= 2\bigl|\{e \in E(A,A)\mid e\notin F\}\bigr|+ \bigl|\{e\in E(A,B)\mid e\notin F\}\bigr|+ \bigl|\{e\in E(A,C)\mid e\notin F\}\bigr|\\[1mm]
	&\qquad\qquad\ge \bigl|\{e\in E(A,B)\mid e\notin F\}\bigr|+ \bigl|\{e\in E(A,C)\mid e\notin F\}\bigr|\\[1mm]
	&\qquad\qquad\ge e(A,B)-2|B_1|-|B_3|+q-2|B_2|-|B_3|\ge e(A,B)+q-2|B|.
\end{align*}
Since we also have $\sum\limits_{x\in A}(d(x)-2)= 2e(A,A)+e(A,C)+e(A,B)-2|A|$, it follows that $q\le 2e(A,A)-2|A|+2|B|+e(A,C)$.
This proves the easy part of Corollary~\ref{2FT}, and thus also of Theorem~\ref{2FTrevised}, the `only if' direction.\hfill$\Box$

\medskip
\noindent
In case of equality, i.e.\ if $q= 2e(A,A)-2|A|+2|B|+e(A,C)$, we have equality in all the above inequalities.
From this we obtain the following.

\clearpage
\begin{thm}\mbox{}\\*
Let $G$ be a graph with a $2$-factor $F$.
Let $A,B\subseteq V(G)$ be disjoint subsets.
Set $C=V(G)\setminus\linebreak(A\cup B)$, and suppose the graph $G[C]$ has $q= 2e(A,A)-2|A|+2|B|+e(A,C)$ connected components joined to $A$ with an odd number of edges.
Construct $B_1,B_2,B_3$ as above.
Then all edges of $E(A,C)$ belong to $F$ except exactly one edge from $q-2|B_2|-|B_3|$ of the components joined to $A$ with an odd number of edges.
Moreover, we have $B=B_1\cup B_2\cup B_3$.
\label{2FT=}
\end{thm}

\subsection{The `if' part of the 2-Factor Theorem}

To prove the `if' part of the $2$-Factor Theorem (Theorem \ref{2FT}) we shall use the theory of alternating chains of Gallai from his 1950 paper \cite{Gallai 1950}.
A main result of that paper was also obtained independently by Belck \cite{Belck 1950}, submitted before Gallai's paper had appeared.
Moreover, we shall use a simple but useful lemma by Belck.
In that sense, the proof of the $2$-Factor Theorem we present is due to Belck.

Let $G$ be a graph, possibly with multiple edges and loops, and let each edge of $G$ be coloured either blue or red.
An \emph{alternating chain} is a walk without repeated edges and with the edges coloured alternating blue and red.
(A single vertex is considered both as a red-blue and a blue-red alternating chain.)
We fix some vertex $p$ and consider alternating chains in~$G$ starting from $p$ with a blue edge.
The vertices of $G$ can then be divided into four groups: 1)~\emph{BR-vertices} can be reached by an alternating chain as described ending in blue and also by an alternating chain ending in red; 2)~\emph{B-vertices} can be reached only by alternating chains ending in blue; 3)~\emph{R-vertices} can be reached only by alternating chains ending in red; and finally 4)~the \emph{unreachable vertices}.

The single vertex $p$ is itself a blue-red alternating chain, so $p$ is either an R-vertex or a BR-vertex.
Belck and Gallai independently proved the following.

\begin{thm}[Belck \cite{Belck 1950}, Gallai \cite{Gallai 1950}]\mbox{}\\*
	Let $C$ be a connected component in the graph induced by the BR-vertices.
	If $p$ does not belong to $C$, then $C$ has exactly one entering edge, either a blue edge from an R-vertex to a vertex in $C$ or a red edge from a B-vertex to a vertex in $C$, such that any alternating chain from $p$ starting with a blue edge to any vertex in $C$ must enter $C$ via the entering edge.
\label{GB}
\end{thm}

\noindent
The proof can be found in Belck \cite[\S\,4]{Belck 1950}, and can be more easily, and with all details, be found in Gallai's paper \cite[\S\,3]{Gallai 1950}.
It was later presented by Berge in `the second book' of graph theory \cite{Berge 1958}.
The idea of the proof is that if there were a second entering edge, then using that edge walking into $C$ and later the first entering edge walking out of $C$ again, the end-vertex of the first entering edge outside $C$ would be a BR-vertex and thus belong to $C$.

We next formulate a lemma by Belck \cite[Zusatz~I,\,1]{Belck 1950} that we use only in some special situations.
(Belck's version applies to $k$-factors in general.)

\begin{lem}[Belck \cite{Belck 1950}]\mbox{}\\*
	Let $G$ be a maximal graph without a $2$-factor, i.e.\ adding any new edge to $G$ (which can be a parallel edge or a loop) will result in a graph with a $2$-factor.
	Let $e_1$ and $e_2$ be two new edges.
	Let $F_1$ be a $2$-factor in $G+e_1$ and $F_2$ a $2$-factor in $G+e_2$.
	Colour the edges in $E(F_1)\setminus E(F_2)$ red and the edges in $E(F_2)\setminus E(F_1)$ blue.
	Then $H=G+e_1+e_2$ contains a closed alternating chain containing $e_1$ and $e_2$.
\label{Belck lemma}
\end{lem}

\noindent
We give an outline of the proof.
Let $e_1=xy$ and start a chain in $H$ from $x$ with the red $E(F_1)\setminus E(F_2)$ edge $xy$.
Then $y$ is incident to a blue $E(F_2)\setminus E(F_1)$ edge $yz$.
Then $z$ is incident to a new red $E(F_1)\setminus E(F_2)$ edge, and so on.
This process can only stop at $x$ with a blue edge from $E(F_2)\setminus E(F_1)$.
The obtained chain must contain $e_2$, since otherwise changing its colours would result in a red $2$-factor of $G$.
This proves the lemma.\hfill$\Box$

\medskip\noindent
A consequence of Lemma \ref{Belck lemma}, which again is a special case of a result for $k$-factors with the case $k = 1$ proved also in \cite{Lovasz 1975}, is the following.

\begin{lem}\mbox{}\\*
Let $G$ be a maximal graph without a $2$-factor.
Let $a,b,c,d$ be vertices in $G$, where $a$ and $b$ are joined by two edges, $b$ and $c$ are also joined by two edges, and $b$ and $d$ are joined by at most one edge.
Then $a$ and $c$ are joined by two edges.
\end{lem}

\noindent
We just give an outline of the proof again.
Assume otherwise, and let $e_1$ and $e_2$ be new edges between $b$ and $d$, and between $a$ and $c$, respectively.
Follow Lemma \ref{Belck lemma} to colour some edges red or blue.
Then, by the lemma, there is a closed alternating chain containing the new red edge $e_1$ and the new blue edge $e_2$.
Hence there is an alternating subchain $S$ from $b$ to either~$a$ or $c$, starting with the red edge $e_1$ and ending with a red edge, say in $a$.
One of the two edges between $a$ and $b$ is not red.
Colour that edge red and change colours on $S$.
The result is a red $2$-factor in $G$, which is a contradiction.\hfill$\Box$

\medskip\noindent
We are now ready to prove the `if' part of the $2$-Factor Theorem \ref{2FTrevised}.
Let $H$ be a graph without a $2$-factor.
Add edges to $H$ to obtain a maximal graph $G$ without a $2$-factor.
Remember that we operate within the class $\mathcal{M}_2$ of graphs with edge-multiplicities at most $2$ and at most one loop at each vertex.
We shall prove that $G$ has an independent set $A\subseteq V(G)$, and a set $B\subseteq V(G)$ disjoint from $A$, such that if we set $C = V(G)\setminus(A\cup B)$, then the number~$q$ of components of $G[C]$ joined to $A$ by an odd number of edges is strictly greater than  $-2|A|+2|B|+e(A,C)$.
Then the same statement holds for $H$, proving the `if' part of Corollary~\ref{2FT} and Theorem \ref{2FTrevised}.

Let $A$ be the set of vertices of $G$ without loops, and let $B$ be the set of vertices with a loop and joined to all other vertices of $G$ by exactly two edges.
Let $C$ be the remaining vertices.
The set $A$ is non-empty, since otherwise each vertex of $G$ would have a loop and all the loops together would form a $2$-factor.
Let $p\in A$ and add a loop $e_1$ at $p$ to $G$.
By maximality of~$G$, the resulting graph has a $2$-factor $F_1$.
Colour the edges in $F_1$ red and all other edges of $G$ blue.
Now consider alternating chains starting in $p$ with blue edges.
(We might alternatively say starting with the red loop at $p$.)
Such an alternating chain cannot end at $p$ with a blue edge, since changing red and blue on it would result in a red $2$-factor of $G$.
We are then in the situation described in the start of this section, and we may use the theory and Theorem~\ref{GB} of Belck and Gallai.

Since $p$ is not reachable by an alternating chain ending with a blue edge, it must be an R-vertex.
And in fact all vertices of $A$ are R-vertices.
Let namely $r\in A$ be a vertex different from $p$.
Then adding a loop $e_2$ at $r$ to $G$ gives a $2$-factor $F_2$.
Using Belck's lemma \ref{Belck lemma} (with~$e_1$ the loop added at $p$ in the previous paragraph, and $e_2$ the new loop at~$r$), we get two edge-disjoint alternating chains in $G$ starting from $p$ with new blue edges and ending in~$r$ with red edges (see Figure~1).
\begin{figure}[htb]\label{fig1}
	\centerline{\begin{tikzpicture}[scale=1,
			vertex/.style={circle, draw, fill=white, inner sep=2pt, minimum size=6pt},
			vertex_sm/.style={circle, draw, fill=white, inner sep=0.5pt, minimum size=2pt},
			edge_F1/.style={line width=2.5pt, red},
			edge_notF1/.style={line width=1.1pt, blue},
			edge_not/.style={line width=1.1pt, black}]

			\node[vertex] (p) at (0,0.5) {};
			\draw[edge_F1] (p) to[out=50,in=130,looseness=18] (p);
			\node at (0,1.55) {$e_1$};

			\node[vertex] (ul1) at (1, 1) {};
			\node[vertex] (ul2) at (2, 1) {};
			\draw[edge_F1] (ul1) -- (ul2);

			\node[vertex] (ml1) at (1, 0) {};
			\node[vertex] (ml2) at (2, 0) {};
			\draw[edge_F1] (ml1) -- (ml2);

			\draw[edge_notF1] (p) -- (ul1);
			\draw[edge_notF1] (p) -- (ml1);

			\node at (-0.3,0) {$p \in A$};
			\node at (-0.3,-0.5) {R-vertex};

			\node[vertex_sm] at (2.5,1) {};
			\node[vertex_sm] at (2.75,1) {};
			\node[vertex_sm] at (3,1) {};
			\node[vertex_sm] at (3.25,1) {};
			\node[vertex_sm] at (3.5,1) {};
			\node[vertex_sm] at (2.5,0) {};
			\node[vertex_sm] at (2.75,0) {};
			\node[vertex_sm] at (3,0) {};
			\node[vertex_sm] at (3.25,0) {};
			\node[vertex_sm] at (3.5,0) {};

			\node[vertex] (r) at (6,0.5) {};
			\draw[edge_not] (r) to[out=50,in=130,looseness=18.5] (r);
			\node at (6,1.55) {$e_2$};

			\node[vertex] (ur1) at (4, 1) {};
			\node[vertex] (ur2) at (5, 1) {};
			\draw[edge_notF1] (ur1) -- (ur2);

			\node[vertex] (mr1) at (4, 0) {};
			\node[vertex] (mr2) at (5, 0) {};
			\draw[edge_notF1] (mr1) -- (mr2);

			\draw[edge_F1] (r) -- (ur2);
			\draw[edge_F1] (r) -- (mr2);

			\node at (6.2,0) {$r \in A$};

			\draw[edge_F1]    (-1,-1.5) -- (0,-1.5) node[right,text=black] {edge in $F_1$};
			\draw[edge_notF1] (3,-1.5) -- (4,-1.5) node[right,text=black] {edge outside $F_1$};

	\end{tikzpicture}}
	\caption*{Figure 1}
\end{figure}
By Belck's lemma, these new blue edges are not in $E(F_1)$, hence they are also blue in the original colouring.
So $r$ is either an R-vertex or a BR-vertex.
If~$r$ is a BR-vertex, then walk from $r$ along the two alternating chains in the direction of $p$.
On both chains we eventually move out of the BR-component containing $r$, but then that component has two entering edges, contradicting the Gallai-Belck theorem \ref{GB}.
Hence $r$ must be an R-vertex, as claimed.

Since all vertices in $A$ are R-vertices, the vertex $p$ is therefore not joined to any $r\in A$ by a blue edge, and not by a red edge either (because of the red loop at $p$).
Hence $p$ has no neighbours in $A$.
Because $p$ was chosen as an arbitrary vertex of $A$, we obtain that $A$ is an independent set.

Take $b\in B$.
Then $b$ is joined to $p$ by two blue edges, one of which may be considered as an alternating chain ending with a blue edge.
Assume there is an alternating chain from $p$ starting with a blue edge and ending at $b$ with a red edge.
If it uses both blue edges between~$p$ and $b$, then change blue and red on it to obtain a $2$-factor of $G$, which is a contradiction.
Otherwise add one of the blue edges to the chain and obtain an alternating chain to the R-vertex $p$ which ends with a blue edge, again a contradiction.
Thus the vertices of $B$ are B-vertices.

The vertices of $C$ have loops, but are not completely joined to everything else.
Take $c\in C$ and let $ct$ be an edge missing in $G$.
As before, by Belck's lemma, there is an alternating chain from~$p$ to~$c$ starting with a blue edge and ending with a red edge.
By definition, there is a loop at~$c$.
If that loop is red, it must be the end-edge of the chain, hence by removing it we get an alternating chain from $p$ to $c$ starting and ending with a blue edge.
If the loop is blue and has been used on the chain, then part of the chain ends in $c$ with a blue edge.
If the loop is blue, but not part of the chain, then add it and get an alternating chain from $p$ to~$c$ starting with a blue edge and ending with a blue edge.
The conclusion is that all vertices of~$C$ are BR-vertices.

The components of $G[C]$ are thus the BR-components.
Since each vertex of $C$ has red degree $2$, the number of red edges from a component to $V(G-C)$ is even.
By Theorem~\ref{GB}, each component has exactly one entering edge.
There are now two possibilities.
If the entering edge is red, then it enters from a B-vertex, all other edges from~$B$ to the component must be blue, and all edges from $A$ to the component must be red and odd in number.
If the entering edge is blue, then it enters from an R-vertex, i.e.\ from a vertex in $A$, all other edges from $A$ to the component must be red and even in number, and all edges from the component to $B$ must be blue since there is only one entering edge.

In conclusion, we get that all components of $G[C]$ have a odd number of edges to $A$.
Suppose that $q_1$ components have a red entering edge and $q_2$ a blue.
Then the number of red edges from $B$ to $A$ equals $2|B|-q_1$.
The number or red edges from $C$ to $A$ equals $e(A,C)-q_2$.
Since the total number of red edges between $A$ and $B\cup C$ is $2|A|-2$, we find
\[2|A|-2= 2|B|-q_1+e(A,C)-q_2= 2|B|+e(A,C)-q,\]
and thus
\[q= -2|A|+2|B|+e(A,C)+2> -2|A|+2|B|+e(A,C).\]
This completes the proof of the `if' part of Corollary~\ref{2FT} and Theorem \ref{2FTrevised}.\hfill$\Box$

\medskip
\noindent
For later use we note the following.
In the maximal graph $G$ without a $2$-factor, the vertex $p\in A$ is joined to $C$ only by blue edges, and thus only by entering edges.
By Theorem~\ref{GB}, there is at most one such edge to each component of $G[C]$, hence $p$ is joined by at most one edge to each component.
Since $p\in A$ was chosen arbitrarily, it follows that each vertex of $A$ is joined to each component of $G[C]$ by at most one edge.

Also, suppose that a vertex $c\in C$ is joined to $p$ by an edge $e_1=pc$, and also by another edge~$e_2$ to $A$.
Since $p$ is joined to each component of $C$ by at most one edge, the edge~$pc$ is single.
If we add a new edge $e_3=pc$ to $G$, rather than adding the loop at $p$, then we obtain a new $2$-factor $F^*$ in $G+e_3$.
Theorem \ref{2FT=} now applies.
(The number of odd components is now one less and the number of edges between $A$ and $C$ one more.)
Thus all edges $e_1,e_2,e_3$ belong to $F^*$, which is impossible, since all three are incident with $x$.
Again, $p\in A$ was chosen arbitrarily, so we can conclude that each vertex of $C$ is joined in $G$ to at most one vertex of~$A$.
The conclusion is that the edges joining vertices in $A$ with the vertices of a component of $G[C]$ form a matching.

Finally, note that in the maximal graph $G$ without a $2$-factor the components of $G[C]$ are all complete in the sense that any two vertices are joined by a double edge and any vertex has a loop.
This follows since adding anything missing still keeps too many components of~$G[C]$ joined to $A$ by an odd number of edges.

\section{Regular graphs without 2-factors}

Now let $G$ be a $(2k+1)$-regular graph without a $2$-factor.
Then the $2$-Factor Theorem \ref{2FTrevised} gives that there exist disjoint subsets $A,B\subseteq V(G)$, $A$ an independent set, such that if we set $C=V(G)\setminus(A\cup B)$, then the number $q$ of connected components of $G[C]$ joined to $A$ by an odd number of edges satisfies
\[q> -2|A|+2|B|+e(A,C)= (2k-1)|A|-e(A,B)+2|B|.\]
Since $q$ and $e(A,C)$ have the same parity, we deduce
\begin{equation}
    q+e(A,B)\ge (2k-1)|A|+2|B|+2.
    \label{first}
\end{equation}
Let $q_1$ denote the number of components of $G[C]$ that are joined to $A$ with exactly one edge and not joined to $B$ at all, and let $q_2$ denote the number of components of $G[C]$ that are joined to $A$ with exactly one edge and to $B$ with at least one edge.
Finally, let~$q_3$ denote the number of components of $G[C]$ that are joined to $A$ with an odd number of at least $3$ edges.
Then of course $q=q_1+q_2+q_3$.
Looking at edges incident with the vertices of $A$, we get
\begin{equation}
   (2k+1)|A|\ge e(A,B)+q_1+q_2+3q_3.
   \label{second}
\end{equation}
Similarly, looking at the edges incident with vertices from $B$, we get
\begin{equation}
    (2k+1)|B|\ge e(A,B)+q_2.
    \label{third}
\end{equation}
Adding $\tfrac12(2k+1)$ times \eqref{first}, $\tfrac12(2k-1)$ times \eqref{second}, and $1$ times \eqref{third} together, we get:
\begin{equation}
    q_1\ge (2k+1)+(2k-2)q_3.
    \label{fourth}
\end{equation}
This proves that $G$ has at least $2k+1$ leaves, which is exactly Theorem \ref{regular}.\hfill$\Box$

\medskip
\noindent
Sylvester, in his correspondence with Petersen in 1889 (see \cite{Konig 1936, Sab 1991}), gave examples of $(2k+1)$-regular graphs without $2$-factors, so-called \emph{primitive graphs}, having exactly $2k+1$ leaves.
Because of Theorem~\ref{JP2}, the primitive graphs are exactly those without $2$-factors.
The above inequalities can be used to describe all the primitive $(2k+1)$-regular graphs having exactly $2k+1$ leaves, because there is then equality in the inequalities.
The case $|A|=1$ in the following theorem gives the Sylvester graphs.

\begin{thm}\mbox{}\\*
Let $G$ be a $(2k+1)$-regular graph, with $k\ge2$, without a $2$-factor and with exactly $2k+1$ leaves.
Then there exists disjoint subsets $A,B\subseteq V(G)$, where $A$ and $B$ are both independent (i.e.\ $A\cup B$ induces a bipartite subgraph of $G$) with $|A|=|B|+1$.
Moreover, $2k+1$ of the connected components of $G-(A\cup B)$ are joined to $A$ with exactly one edge and not to $B$, while all other components are joined to $A$ with exactly one edge and to $B$ with exactly one edge.
\label{primitive}
\end{thm}

\noindent
We may construct these graph as follows.
Start from a bipartite graph $H$ with $A$ and $B$ the sets of vertices in its two sides, where $|A|=|B|+1$ and with all vertices of $B$ of degree $2k+1$ and all vertices of $A$ of degree at most $2k+1$.
Add $2k+1$ disjoint graphs to $H$, each having all vertices of degree $2k+1$ except one vertex of degree $2k$ and join each to $A$ by one edge to get a $(2k+1)$-regular graph with $2k+1$ leaves.
(To get exactly $2k+1$ leaves, we should use $2k+1$ disjoint graphs obtained as follows.
Take a connected $(2k+1)$-regular graph with exactly two leaves and remove a bridge.
A connected component in the remaining graph is suitable as one of the $2k+1$ graphs.)
Moreover, for a number of edges $xy$ of the bipartite graph~$H$, replace the edge $xy$ by a  $(2k+1)$-regular bridgeless graph, removing an edge from it, and joining its two endvertices to $x$ and to $y$ to obtain a $(2k+1)$-regular graph.

In the case of $3$-regular graphs, i.e.\ if $k=1$, we may have $q_3>0$ and still have equality in \eqref{fourth}, hence we have to add the following to the above description: Finally in addition, for some vertices $x \in B$  joined to three vertices of $A$, replace the vertex $x$ by a $3$-regular bridgeless graph, removing a vertex from it and joining its three neighbours to the three neighbours of~$x$ in $A$.

\section{The maximal graphs without 2-factors}

Remember that we operate within the class $\mathcal{M}_2$ of graphs with edges of multiplicity at most~$2$ and with at most one loop at each vertex.
Let $G$ be a graph in $\mathcal{M}_2$ and assume that $G$ does not have a $2$-factor, but that $G$ is maximal with respect to this property, i.e.\ if a new edge or a new loop is added and we stay within the class $\mathcal{M}_2$, then the resulting graph has a $2$-factor.

From our proof of the `if' part of the $2$-Factor Theorem, we immediately get the following partial characterisation, where the only missing part is a precise description of how the vertices of $C$ are joined to the vertices of $A$.

\begin{thm}\mbox{}\\*
	Let $\mathcal{M}_2$ be the class of all graphs with all multiplicities at most $2$ and with each vertex having at most one loop, and let $G$ be maximal without a $2$-factor within the class $\mathcal{M}_2$.
	Let the set~$A$ contain all vertices of $G$ without loops, the set $B$ contain all vertices in $G$ with loops and joined to all other vertices in $G$ by two edges, and set $C=V(G)\setminus(A\cup B)$.
	Suppose $G[C]$ has~$q$ connected components $\mathcal{C}=\{C_1,\ldots,C_q\}$.
	Then we have the following.
	\vspace{-2mm}
	\begin{itemize}
		\addtolength{\itemsep}{-2mm}
		\item The set $A$ is independent.
		\item The components in $\mathcal{C}$ are all complete (with a loop at every vertex and two edges between any two vertices in the same component).
		\item Each component $C_i$ is joined to $A$ by an odd matching (of size $e(A,V(C_i))$).
		\item $2|A|+q=2|B|+e(A,C)+2$, which is equivalent to $|A|= |B|+1+\sum\limits_{i=1}^q\tfrac12\bigl(e(A,V(C_i))-1)$.
	\end{itemize}
\label{2FTmax}
\end{thm}

\noindent
The case $q=0$ is possible; in this case $|A|=|B|+1$.
The resulting graph is maximal.
(Remember that $A$ is an independent set and each vertex of $B$ has a loop and is joined to all other vertices by two edges.)
For $q=1$ and $q=2$, the conditions described above are necessary and sufficient for $G$ to be maximal without a $2$-factor.
For $q\ge3$ we need to add extra conditions; for example for $q=3$:
\vspace{-2mm}
\begin{itemize}
\item \textit{No three components of $G[C]$ are all joined to the same two different vertices $x,y\in A$.}
\end{itemize}
Assume namely that this was the case.
Then adding a loop at $x$ gives a graph $G^*$ having a $2$-factor $F$ and satisfying Theorem \ref{2FT=}.
The three edges from the components to $x$ are not in~$F$, because the loop at $x$ is.
But then the three edges to $y$ are all in $F$, a contradiction.

In what follows we shall generalise this property to a general necessary and sufficient condition for $G\in\mathcal{M}_2$ to be maximal without a $2$-factor, as exhibited in the final condition of Theorem~\ref{final}.
The above is the case with $A'=\{x,y\}$ and $C'=\{C_1,C_2,C_3\}$.

\subsection{The necessary and sufficient conditions}

Let again $G$ be a graph with all multiplicities at most $2$ and with each vertex having at most one loop.
Let us assume that $G$ does not have a $2$-factor, and define $A,B,C,q,\mathcal{C}$ as in Theorem \ref{2FTmax}.
For $i=1,\ldots,q$, set $t_i=e(A,C_i)$.
Also assume that the four conditions in the theorem are satisfied.

In this subsection we shall determine some extra properties that allow us establish necessary and sufficient conditions for $G$ to be maximal within $\mathcal{M}_2$ without a $2$-factor.
There are four possible types of additions of a new edge $e$ to $G$: 1)~a loop at a vertex  $x\in A$; 2)~an edge between two different vertices $x,y\in A$; 3)~an edge between a vertex $x\in A$ and a vertex $y$ in a component of $\mathcal{C}$; and, finally, 4)~an edge between two vertices~$x$ and $y$ in two different components of $\mathcal{C}$.

In all four situations, the graphs $G^*= G+e$, with the subsets $A$ and $B$ defined as above, satisfy the condition of Theorem \ref{2FT=} that $q=2e(A,A)-2|A|+2|B|+e(A,C)$.
We shall now in each case find a necessary and sufficient condition for $G^*$ to contain a $2$-factor.

The conditions will be given in terms of properties of a bipartite graph $H$.
The vertices in the two sides of $H$ are $A$ and $\mathcal{C}$, respectively.
Moreover, $x\in A$ and $C_i\in\mathcal{C}$ are joined by an edge in $H$ if and only if $x\in A$ is joined to $C_i$ in $G$.
The bipartite graph $H$ is uniquely determined from $G$, and vice versa (up to isomorphism).
There is a one-to-one correspondence between the edges of $H$ and the edges between $A$ and $G[C]$ in $G$, where $(x,C_i)\in E(H)$ corresponds to an edge $(x,y)\in E(G)$ with $y\in V(C_i)$.
From Theorem~\ref{2FTmax} we have that the vertex $C_i$ has odd degree $t_i$ in $H$.

\medskip
\noindent
\textbf{Case 1: $G^*$ is the graph obtained from $G$ by adding a loop to a vertex $x \in A$.}\\*
Suppose there is a $2$-factor $F^*$ in $G^*$ containing the loop at $x$.
The edges from $F^*$ joining vertices from $A$ to vertices of $C$ correspond to a set of edges $F$ of the bipartite graph~$H$.

By Theorem \ref{2FT=}, $F$ satisfies the following.\\
(1.1) The vertex $x\in A$ is not incident with an edge from $F$.\\
(1.2) Each vertex $y\in A\setminus\{x\}$ is incident with at most two edges from $F$.\\
(1.3) For each $C_i\in \mathcal{C}$, all $t_i$ edges incident with $C_i$, except perhaps one, are in $F$.

If all edges in $H$ incident with $C_i$ are in $F$, then we remove one of the edges for each such~$C_i$ from $F$ and obtain a subset $F'$ of $E(H)$ satisfying the following.\\
(1.1') The vertex $x\in A$ is not incident with an edge from $F'$.\\
(1.2') Each vertex $y\in A\setminus\{x\}$ is incident with at most two edges from $F'$.\\
(1.3') For each $C_i\in \mathcal{C}$, all $t_i$ edges incident with $C_i$, except exactly one, are in $F'$.

We have shown above that if a loop added to the vertex $x\in A$ in $G$ gives rise to a $2$-factor in the resulting graph $G^*$, then $E(H)$ has a subset $F'$ satisfying (1.1'), (1.2') and (1.3').
We shall now show that the converse is true as well.

So assume that $E(H)$ has a subset $F'$ satisfying (1.1'), (1,2') and (1,3').
Then $|F'|=\sum_{i=1}^q(t_i-1)$.
For $i=0,1,2$, let $A_i$ denote the set of vertices in $A$ incident with $i$ edges from $F'$.
Then $|F'|=|A_1|+2|A_2|=\sum_{i=1}^q(t_i-1)$.
Using Theorem \ref{2FTmax}, we deduce that $|B|=|A_0|+\tfrac12|A_1|-1$.
A $2$-factor $F^*$ in the graph $G^*$ can now be obtained as follows.
$F^*$ consists of the loop at $x$, the edges in $G$ corresponding to the edges $F'$ in $H$, for each vertex of $A_0\setminus\{x\}$ a double-edge to a vertex in~$B$, and for each vertex of  $A_1$ one edge to a vertex in~$B$.
Because of the size of $B$, this can be done in such a way that each vertex of $B$ becomes incident to exactly two edges of $F^*$.
For each component $C_i\in\mathcal{C}$, there are now $t_i-1$ matching edges from $F^*$ coming in from $A$.
We finally add the endvertices in $C_i$ two and two by an edge of $C_i$ to $F^*$, and add the loop at all other vertices of $C_i$ to $F^*$.
The result is a $2$-factor~$F^*$ of $G^*$.

In conclusion, the existence of a subset $F'$ of edges of the graph $H$ satisfying (1.1'), (1.2'), (1.3') is equivalent to the graph $G^*$ having a $2$-factor, where $G^*$ is obtained from $G$ by adding a loop at the vertex $x \in A$.

\medskip
\noindent
\textbf{Case 2: $G^*$ is the graph obtained from $G$ by adding an edge between two different vertices $x,y\in A$.}\\*
By arguments similar to those in Case 1, the graph $G^*$ has a $2$-factor if and only if the graph~$H$ has a subset $F'$ of edges satisfying the following.\\
(2.1') The vertices $x,y\in A$ are each incident with at most one edge from $F'$.\\
(2.2') Each vertex $z\in A\setminus\{x,y\}$ is incident with at most two edges from $F'$.\\
(2.3') For each $C_i\in \mathcal{C}$, all $t_i$ edges incident with $C_i$, except exactly one, are in $F'$.

\medskip
\noindent
\textbf{Case 3: $G^*$ is the graph obtained from $G$ by adding an edge between a vertex $x \in A$ and a vertex $y \in V(C_i)$, for some component $C_i\in\mathcal{C}$.}\\*
By arguments similar to those in Case 1, the graph $G^*$ has a $2$-factor if and only if the graph~$H$ has a subset $F'$ of edges satisfying the following.\\
(3.1') The vertex $x\in A$ is incident with at most one edge from $F'$.\\
(3.2') Each vertex $z\in A\setminus\{x\}$ is incident with at most two edges from $F'$.\\
(3.3') For each $C_j\in \mathcal{C}$ with $j\neq i$, all $t_j$ edges incident with $C_j$, except exactly one, are in~$F'$.\\
(3.4') All $t_i$ edges incident with $C_i$ are in $F'$.

\medskip
\noindent
\textbf{Case 4: $G^*$ is the graph obtained from $G$ by adding an edge between a vertex $x\in V(C_i)$ and a vertex $y\in V(C_j)$, for two different components $C_i,C_j\in\mathcal{C}$.}\\*
By arguments similar to that in Case 1, the graph $G^*$ has a $2$-factor if and only if the graph~$H$ has a subset $F'$ of edges satisfying the following.\\
(4.1') Each vertex $x\in A$ is incident with at most two edges from $F'$.\\
(4.2') For each $C_k\in \mathcal{C}$ with $k\notin\{i,j\}$, all $t_k$ edges incident with $C_k$, except exactly one, are in $F'$.\\
(4.3') All $t_i$ edges incident with $C_i$ are in $F'$.\\
(4.4') All $t_j$ edges incident with $C_j$ are in $F'$.

\medskip
As we shall see in the next subsection, the four sets of conditions are, not surprisingly, closely related.
In fact, the first set of conditions for all $x\in A$ imply the second, third and fourth sets of conditions.

\subsection{Reformulation of the conditions; putting it all together}

Let $G$ be a maximal graph without a $2$-factor within the class $\mathcal{M}_2$ satisfying Theorem \ref{2FTmax} and the extra conditions of the previous subsection.
Let again $H$ be the bipartite graph with vertex sets $A$ and $\mathcal{C}$, where $\mathcal{C}$ has $q$ vertices $C_1,\ldots,C_q$ of degrees $t_1,\ldots,t_q$, all odd integers.
We direct all edges from $A$ to $\mathcal{C}$ and give them capacity $1$, add a source and join it by a directed edge  of capacity $2$ to each vertex of $A$, except to one vertex $x\in A$, and add a sink with a directed edge of capacity $t_i-1$ from each $C_i\in \mathcal{C}$.
In this way we obtain from the bipartite graph $H$ a flow network $N$ with integer capacities.
The existence of a subset $F'$ of the edges of $H$ satisfying the conditions (1.1'), (1.2') and (1.3') is then equivalent to the existence of a flow of size $\sum_{i=1}^q (t_i-1)$ in $N$.
(If there is a flow of that size there is also an integer flow of that size, and the existence of such an integer flow may be translated directly into the existence of $F'$.)
By the Max Flow Min Cut Theorem \cite{FF 1962} the existence of  a flow of this size is equivalent to any cut of $N$ having capacity at least $\sum_{i=1}^q (t_i-1)$.

A cut in $N$ may be thought of as  subsets $A'\subseteq A$ and $\mathcal{C}'\subseteq \mathcal{C}$, where $A'$ and $\mathcal{C}'$ are the parts of $A$ and $\mathcal{C}$ belonging to the same side of the cut as the sink.
We may assume that $x$ belongs to $A'$, since moving it there from $A\setminus A'$ can only give the same or smaller capacity.
The capacity of the cut is then
\[\textstyle 2|A'|-2+e(A\setminus A',\mathcal{C}')+\sum\limits_{C_i\in \mathcal{C}\setminus \mathcal{C}'}(t_i-1).\]
So the condition for the existence of $F'$ is
\[\textstyle 2|A'|-2+e(A\setminus A',\mathcal{C}')+\sum\limits_{C_i\in \mathcal{C}\setminus \mathcal{C}'}(t_i-1)\ge \sum\limits_{C_i\in \mathcal{C}}(t_i-1),\]
which is equivalent to
\[\textstyle 2|A'|\ge 2+\sum\limits_{C_i\in \mathcal{C}'}(t_i'-1),\]
where $t_i'$ is equal to the number of neighbours from $A'$ the vertex $C_i\in \mathcal{C}$ has in $H$.
Since the last inequality holds for all $x\in A'$, we can write the condition in the following way:
\begin{itemize}
	\item[(4.1)]
	for all $A'\subseteq A$, $A'\ne\varnothing$, and all $\mathcal{C}'\subseteq \mathcal{C}$ we have:\quad $2|A'|\ge 2+\sum\limits_{C_i\in \mathcal{C}'}(t_i'-1)$.
\end{itemize}

The condition in (4.1) is equivalent to the existence of $F'$ in Case 1 above.
In a similar way, conditions equivalent to the existence of $F'$ in Cases 2, 3 and 4 can be obtained:
\begin{itemize}
	\addtolength{\itemsep}{-2mm}
	\item[(4.2)]
	for all $A'\subseteq A$, $|A'|\ge2$, and all $\mathcal{C}'\subseteq \mathcal{C}$ we have:\quad $2|A'|\ge 2+\sum\limits_{C_i\in \mathcal{C}'}(t_i'-1)$;
	\item[(4.3)]
	for all $A'\subseteq A$, $A'\ne\varnothing$, and all $\mathcal{C}'\subseteq \mathcal{C}$, $\mathcal{C}'\ne\varnothing$, we have:\quad $2|A'|\ge 2+\sum\limits_{C_i\in \mathcal{C}'}(t_i'-1)$;
	\item[(4.4)]
	for all $A'\subseteq A$ and all $\mathcal{C}'\subseteq \mathcal{C}$, $|\mathcal{C}'|\ge 2$, we have:\quad $2|A'|\ge 2+\sum\limits_{C_i\in \mathcal{C}'}(t_i'-1)$.
\end{itemize}

All these three conditions are implied by (4.1), completing the proof of Theorem~\ref{final}.\hfill$\Box$

\section{Conclusion}

Cycles in graphs and $2$-factors are important structures in graph theory, so even if partial characterisations of the maximal graphs without $2$-factors have existed in the mathematical literature since 1950, it seems worthwhile to obtain a complete characterisation, as has been one of the goals of this paper.
We have also pointed out the important role maximal graphs have played in the development of factorisation theory, and how the theorem of Belck and Gallai (Theorem~\ref{GB}) may be used to prove the $2$-Factor Theorem, and indeed the general $k$-factor theorem.

\medskip
For general $k$-factors, $k\ge3$, it is possible to carry out a similar analysis as done above, but the outcome will be much more complicated.
A component $C_i$ of $G[C]$ is said to be odd if the number $k|V(C_i)|+e(A,V(C_i))$ is odd.
The connections between the set $A$ and the components $C_i$ are no longer restricted to be matchings.

\medskip
Any simple graph $H$ that is maximal without a $2$-factor within the class of simple graphs is equal to the simple restriction of a graph $G$ as described in Theorem~\ref{final}, i.e.\ $H$ is equal to~$G$ with each multiple edge reduced to a single edge and all loops removed.
But not all such reduced versions of graphs satisfying Theorem~\ref{final} are maximal without a $2$-factor within the class of simple graphs.
When $|V(G)|\le4$ there is exactly one exception: $G\cong K_{1,3}$, with $|A|=1$, $|B|=0$ and $q=3$.
The exceptions for $|V(G)|\ge5$ are exactly the following three (overlapping) types of graphs $G$:
\vspace{-2mm}
\begin{itemize}
	\addtolength{\itemsep}{-2mm}
  \item there is a component $C_i\in\mathcal{C}$ such that $e(A,V(C_i))=1$ and $|V(C_i)|\in\{2,3\}$;
  \item $|B|=1$, and there is an $x\in A$ that is not joined to any component in $\mathcal{C}$;
  \item there is an $x\in A$ joined to at least two components $C_i,C_j\in\mathcal{C}$, $j\ne i$, with $|V(C_i)|=1$.
\end{itemize}
\vspace{-2mm}
Proving that this covers all exceptions requires some tedious but fairly straightforward case checking.

\medskip
In some cases the $2$-factor obtained in a maximal graph after adding a new edge may be chosen as a Hamiltonian cycle.
It might be interesting to characterize those simple graphs, maximal without a $2$-factor, that are also maximal without a Hamiltonian cycle.

\section{Hans-Boris Alexander Belck, 1929--2007}

At the conference \emph{Combinatorics in Cambridge} in August 2003, Roland H\"aggkvist gave a lecture \emph{Factors Galore}~\cite{Roland 2001}.
In his abstract he mentioned and asked for data of Hans-Boris Belck who in 1950 published his one and only mathematical paper, a 25-page opus in \emph{Journal f\"ur die Reine und Angewandt Mathematic} (often called \emph{Crelle's Journal}).
In this paper Belck not only gave the first general $k$-factor criterium, including the first purely graph-theoretical proof of Tutte's $1$-Factor Theorem, but also settled the following question completely: For what values of $(l,r,k)$ is an $l$-edge-connected $r$-regular graph (of even order if $k$ is odd) guaranteed to have a $k$-factor?

The second author took an interest in H\"aggkvist's question, and a few years later, with the help of the internet, he stumbled upon a Belck patent from Brazil.
Searching then for Belck families in Brazil, he got in contact with Hans-Boris Belck's two sons, one living in Sao Paulo in Brazil, the other in Germany.
They provided details from the life of their father and explained, at least partially, why the paper \cite{Belck 1950} was his only mathematical contribution.

We now know the following.
\vspace{-2mm}
\begin{itemize}
	\addtolength{\itemsep}{-2mm}
  \item Hans-Boris Alexander Belck was born January 19, 1929, in Apolda, Th\"uringen, Germany.
  \item In 1931 the family moved to Bromberg (now Bydgoszcz in Poland), where they stayed almost until the end of World War II.
  \item He finished high school at age 16 and was admitted to the University of Frankfurt.
  Later he had a 1-year stay at the University of Bern.
  \item In December 1948 he obtained the degree of Diplom-Mathematiker from the University of Frankfurt; and in June 1949 (aged 20) the degree of Doktor der Naturwissenschaften from the same university (thesis title: \emph{Regul\"are Faktoren von Graphen}).
  Supervisors were Ruth Moufang and Wolfgang Franz.
  \item His paper \cite{Belck 1950} was submitted in November 1949, containing the theory of alternating chains and the $k$-factor theorem.
  \item Belck emigrated to the United States in 1950.
  He worked for the  General Electric Company in Massachusetts.
  \item He became a member of AMS as Mr.\ Hans Boris Belck (since as a 21-year-old his doctorate was not recognized).
  \item In 1954 he obtained a Ph.D. in physics from Rensselaer Polytechnic Institute (thesis title: \emph{The Application of a Magnetic Tape Recorder in Analog Computing}).
  He submitted a US patent application in 1953, which was granted in 1960: \emph{System for recording and reproducing signal waves}.
  \item Belck moved to Brazil in 1956.
  He married in 1957 and had three children (Andreas 1960, Monica 1962, and Alexander 1964).
  \item In 1964 he founded Amelco S A Industria Electronika.
  He also obtained several patents in Brazil.
  \item Belck was involved in the development of the Plano Cruzado, an economical plan to fight inflation, which was adopted by the Brazilian government in 1986.
  \item He visited Germany and the University of Frankfurt in 1987, where he had a hearty rendezvous with Wolfgang Franz.
  \item Belck passed away September 29,  2007, aged 78, in Sao Paulo, Brazil.
\end{itemize}

\subsubsection*{Acknowledgement}
We like to thank two anonymous referees for careful reading and for suggestions that greatly improved the presentation in our paper.

Thanks are also due to Andreas and Alexander Belck for information about their father Hans-Boris Alexander Belck.

\end{document}